\documentclass{revtex4}
\usepackage{graphicx}
\usepackage{amssymb}
\usepackage{amsthm}
\usepackage{amsmath}
\begin{document}
\title{Nonlinear stationary differential equations with fractional order $1<\alpha<2$ on metric star graphs}
\author{K.K.Sabirov}

\affiliation{Institute of Mathematics named after V.I.Romanovsky of the Academy of Sciences of the Republic of Uzbekistan, University street 9, Tashkent, 100174, Uzbekistan}

\begin{abstract}
  In this paper we consider nonlinear stationary fractional-in-space differential equations with order $1<\alpha<2$ on the metric star graph with three finite bonds. At the branched point of the star graph we put the weight continuity and the generalized Kirchhoff rule. It is found the exact solutions of nonlinear stationary fractional equations on the star graph. These can be extended to the star graphs with any number of bonds.
\end{abstract}
\maketitle

\section{Introduction} 

In the past few decades of years it has attracted much attentions to the fractional differential equations due to applications in numerous seemingly diverse and widespread fields of science and engineering \cite{Hilfer}-\cite{Umarov}. It is in constant development, which is fed by ideas and results of various directions in mathematics. One of these directions is nonlinear ordinary and partial differential equations considered on the branched structures so-called metric graphs.

The application of differential equations on metric graphs in quantum theory is described in detail in the works \cite{Smilansky}-\cite{Exner}. In \cite{Zarif}, conservation laws, such as energy and norm were studied for deriving boundary conditions at the vertices of metric graphs and traveling soliton solution was considered. Static soliton solutions of stationary nonlinear Schr\"odinger equations with the boundary condition at the vertices of metric graphs were studied in the Refs. \cite{Adami2011}-\cite{Karim}. Static nonlinear waves in networks described by a nonlinear Schr\"odinger equation with point-like nonlinearities on metric graphs were studied, explicit solutions fulﬁlling vertex boundary conditions were obtained and spontaneous symmetry breaking caused by bifurcations was found in the Ref. \cite{Karim2022}.

In this paper we study two types of nonlinear fractional-in-space stationary equations on the metric star graph with three finite bonds. In the section III we give formulation of the problem described by nonlinear fractional stationary differential equation with order $1<\alpha<2$ and we put the boundary conditions as the weight continuity and the generalized Kirchhoff rule at the branched vertex. In the section IV we obtain the static solution on the metric star graph fulfilling after gluing boundary conditions at the vertices. In the section V we consider another nonlinear stationary fractional-in-space equation with order $1<\alpha<2$ on the star graph.

\section{Preliminaries}

{\bf Definition 1:} Let $\Omega=[a;b]$ be finite interval on real axis $R$. The Riemann-Lioville fractional integrals denoted by $I_{a,x}^\alpha$ and $I_{x,b}^\alpha$ of order $\alpha\in C\, (Re(\alpha)>0)$ are defined as
\begin{equation}
\left(I_{a,x}^\alpha f\right)(x):=\frac{1}{\Gamma(\alpha)}\underset{a}{\overset{x}{\int}}\frac{f(t)dt}{(x-t)^{1-\alpha}},\,(x>a;\,Re(\alpha>0))\label{rli1}
\end{equation}
and
\begin{equation}
\left(I_{x,b}^\alpha f\right)(x):=\frac{1}{\Gamma(\alpha)}\underset{x}{\overset{b}{\int}}\frac{f(t)dt}{(t-x)^{1-\alpha}},\,(x<b;\,Re(\alpha>0))\label{rli2}
\end{equation}
respectively, where $\Gamma(\alpha)$ is the Gamma function, integrals (\ref{rli1}) and (\ref{rli2}) are called the left-sided and the right-sided fractional integrals.

{\bf Definition 2} The Riemann-Lioville derivatives denoted $D_{a,x}^\alpha y$ and $D_{x,b}^\alpha y$ of order $\alpha\in {\bf C}$ $(Re(\alpha)\geq 0)$ are defined by
\begin{equation}
\left(D_{a,x}^\alpha y\right)(x):=\frac{1}{\Gamma(1-\alpha)}\left(\frac{d}{dx}\right)^n\underset{a}{\overset{x}{\int}}\frac{y(t)dt}{(x-t)^{\alpha-n+1}},\,(n=\left[Re(\alpha)\right]+1,\,x>a),\label{rlfd1}
\end{equation}
and
\begin{equation}
\left(D_{x,b}^\alpha y\right)(x):=\frac{1}{\Gamma(1-\alpha)}\left(-\frac{d}{dx}\right)^n\underset{x}{\overset{b}{\int}}\frac{y(t)dt}{(t-x)^{\alpha-n+1}},\,(n=\left[Re(\alpha)\right]+1,\,x<b),\label{rlfd2}
\end{equation}
respectively, where $[Re(\alpha)]$ means integer part of $Re(\alpha)$. In particular, $0<\alpha<1$, then
\begin{equation}
\left(D_{a,x}^\alpha y\right)(x)=\frac{1}{\Gamma(1-\alpha)}\frac{d}{dx}\underset{a}{\overset{x}{\int}}\frac{y(t)dt}{(x-t)^\alpha},\,(x>a),\label{rlfdf1}
\end{equation}
and
\begin{equation}
\left(D_{x,b}^\alpha y\right)(x)=-\frac{1}{\Gamma(1-\alpha)}\frac{d}{dx}\underset{x}{\overset{b}{\int}}\frac{y(t)dt}{(t-x)^\alpha},\,(x<b)\label{rlfdf2}
\end{equation}

\section{Formulation of the problem.}

We consider the star graph $G$ with tree bonds $e_j,\,j=1,2,3$, for which a coordinate $x_j$ is assigned. Choosing the origin of coordinates at the vertex, $L$, for bond $e_j$ we put $x_j\in[0;L_j]$. We use the shorthand notation $y_{j}(x)$ for $y_j(x_j)$ where $x$ is the coordinate on the bond $j$ to which the component $y_j$ refers. 

\begin{figure}[ht!]
\centering
\includegraphics[scale=0.3]{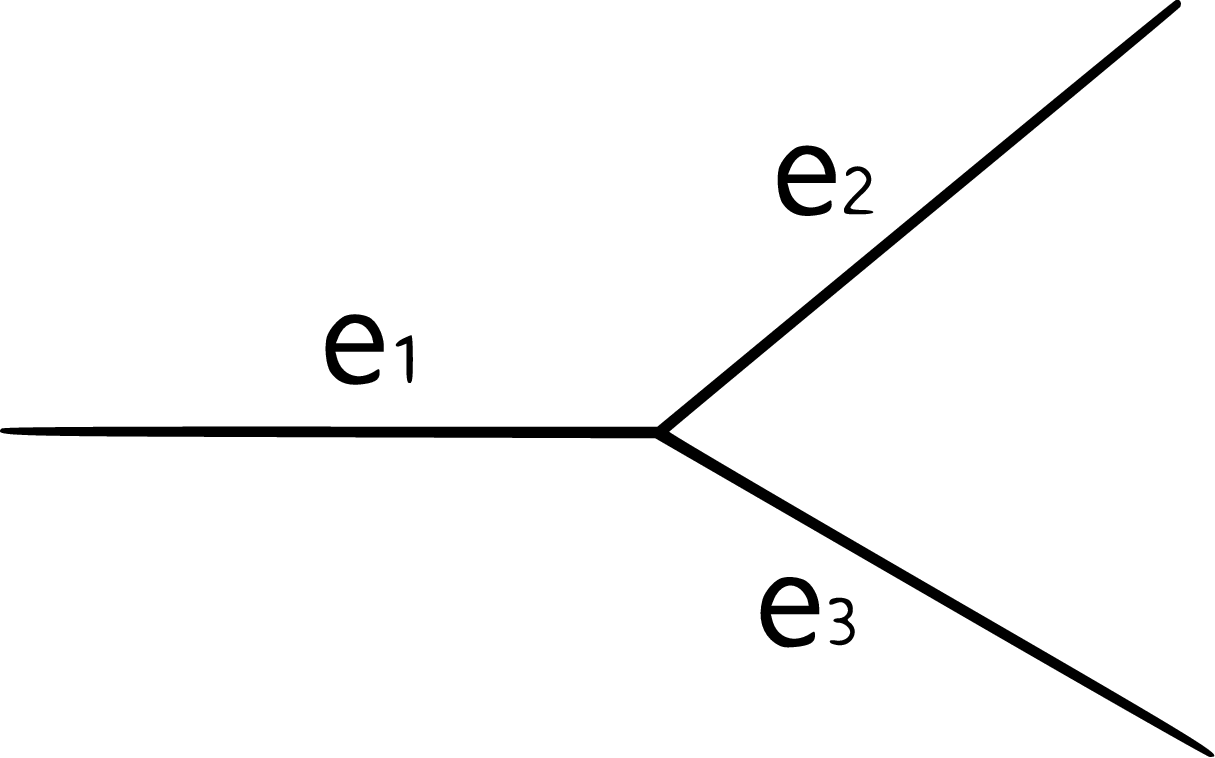}
\caption{The metric star graph} \label{pic1}
\end{figure}

We consider the following the nonlinear stationary differential equation with fractional order $1<\alpha<2$ on the each bond 

\begin{eqnarray}
\left(D_{0,x}^\alpha y_j\right)(x)=\lambda_j x^{\beta_j}\left[y_j(x)\right]^{m_j},\,0<x<L_j,\,m_j\not=1,\,j=1,2,3\label{nsfeq}
\end{eqnarray}
with boundary conditions at the branched point
\begin{eqnarray}
\lambda_1^{\frac{1}{m_1-1}}y_1(x)|_{x=L_1}=\lambda_2^{\frac{1}{m_2-1}}y_2(x)|_{x=L_2}=\lambda_3^{\frac{1}{m_3-1}}y_3(x)|_{x=L_3},\label{bc1}\\
\lambda_1^{\frac{m_1}{m_1-1}}\left(D_{0,x}^{\alpha-1}y_1\right)(L_1-)=\lambda_2^{\frac{m_2}{m_2-1}}\left(D_{0,x}^{\alpha-1}y_2\right)(L_2-)+\lambda_3^{\frac{m_3}{m_3-1}}\left(D_{0,x}^{\alpha-1}y_3\right)(L_3-),\label{bc2}
\end{eqnarray}
and the end of bonds
\begin{eqnarray}
\left(I_{0,x}^{2-\alpha}y_j\right)(+0)=0,\,j=1,2,3,\label{ebc1}
\end{eqnarray}
where $y_j(x)\in C_{\gamma_j^*-\alpha}[0,L_j]$ for $0 <\gamma_j^*-\alpha< 1,\,\gamma_j^*=\frac{\beta_j+m_j\alpha}{1-m_j},\,\lambda\not=0$.

{\bf Theorem 1:} The problem (\ref{nsfeq})-(\ref{ebc1}) has the exact solution.

\section{Proof of the Theorem 1.}
We look for a solution of equation (\ref{nsfeq}) in the form
\begin{eqnarray}
y_j(x)=A_jx^{\gamma_j-1},\,j=1,2,3,\label{esol1}
\end{eqnarray}
where $A_j$ is constant. After directly putting (\ref{esol1}) into (\ref{nsfeq}) the exact solution of equation (\ref{nsfeq}) is the following
\begin{eqnarray}
y_j(x)=\left[\frac{\Gamma\left(\frac{\beta_j+\alpha}{1-m_j}+1\right)}{\lambda\Gamma\left(\frac{\beta_j+m_j\alpha}{1-m_j}+1\right)}\right]^{\frac{1}{m_j-1}}x^{\frac{\beta_j+\alpha}{1-m_j}}.\label{esol2}
\end{eqnarray}

Fulfilling the vertex boundary conditions (\ref{bc1})-(\ref{bc2}) we have the following system of transcendental equations with respect to $\lambda_j,\,j=1,2,3$
\begin{eqnarray}
&\left[\frac{\Gamma\left(\frac{\beta_1+\alpha}{1-m_1}+1\right)}{\Gamma\left(\frac{\beta_1+m_1\alpha}{1-m_1}+1\right)}\right]^{\frac{1}{m_1-1}}L_1^{\frac{\beta_1+\alpha}{1-m_1}}=\left[\frac{\Gamma\left(\frac{\beta_2+\alpha}{1-m_2}+1\right)}{\Gamma\left(\frac{\beta_2+m_2\alpha}{1-m_2}+1\right)}\right]^{\frac{1}{m_2-1}}L_2^{\frac{\beta_2+\alpha}{1-m_2}}=\left[\frac{\Gamma\left(\frac{\beta_3+\alpha}{1-m_3}+1\right)}{\Gamma\left(\frac{\beta_3+m_3\alpha}{1-m_3}+1\right)}\right]^{\frac{1}{m_3-1}}L_3^{\frac{\beta_3+\alpha}{1-m_3}},\label{streq1}\\
&\frac{\lambda_1}{\left(\frac{\beta_1+m_1\alpha}{1-m_1}+1\right)}\left[\frac{\Gamma\left(\frac{\beta_1+\alpha}{1-m_1}+1\right)}{\Gamma\left(\frac{\beta_1+m_1\alpha}{1-m_1}+1\right)}\right]^{\frac{1}{m_1-1}+1}L_1^{\frac{\beta_1+\alpha}{1-m_1}+1}=\frac{\lambda_2}{\left(\frac{\beta_2+m_2\alpha}{1-m_2}+1\right)}\left[\frac{\Gamma\left(\frac{\beta_2+\alpha}{1-m_2}+1\right)}{\Gamma\left(\frac{\beta_2+m_2\alpha}{1-m_2}+1\right)}\right]^{\frac{1}{m_2-1}+1}L_2^{\frac{\beta_2+\alpha}{1-m_2}+1}+\nonumber\\
&\frac{\lambda_3}{\left(\frac{\beta_3+m_3\alpha}{1-m_3}+1\right)}\left[\frac{\Gamma\left(\frac{\beta_3+\alpha}{1-m_3}+1\right)}{\Gamma\left(\frac{\beta_3+m_3\alpha}{1-m_3}+1\right)}\right]^{\frac{1}{m_3-1}+1}L_3^{\frac{\beta_3+\alpha}{1-m_3}+1}.\label{streq2}
\end{eqnarray}
In the particular case when $m_j=m,\,\beta_j=\beta,\,L_j=L,\,j=1,2,3$, we have $\lambda_1=\lambda_2+\lambda_3$, therefore the system of transcendental equations (\ref{streq1})-(\ref{streq2}) has a solution. Now it is clear that the exact solution (\ref{esol2}) fulfills not only the boundary conditions (\ref{ebc1}) and satisfies the following conditions
\begin{eqnarray}
\left(D_{0,x}^{\alpha-1}y_j\right)(+0)=0,\,j=1,2,3.\label{ebc2}
\end{eqnarray}
Therefore the considered problem has the exact solution in the form (\ref{esol2}). Theorem is proved.
\section{Another nonlinear stationary fractional equation on the star graph}

We consider the following the nonlinear stationary differential equation with fractional order $1<\alpha<2$ on the each bond 
\begin{eqnarray}
\left(D_{0,x}^\alpha y_j\right)(x)=\lambda_j x^{\beta_j}\left[y_j(x)\right]^{m_j}+b_jx^{\nu_j},\,0<x<L_j,\,m_j>0,\,m_j\not=1,\,j=1,2,3\label{nsfeq2}
\end{eqnarray}
with weight continuity boundary condition (\ref{bc1}) and (\ref{bc2}) at the branched point and the conditions (\ref{ebc1}) at the end of bonds, where $y_j(x)\in C_{\gamma_j^*-\alpha}[0,L_j]$ for $0 <\gamma_j^*-\alpha< 1,\,\gamma_j^*=\frac{\beta_j+m_j\alpha}{1-m_j},\,\lambda\not=0$.

We look for a solution of equation (\ref{nsfeq}) in the form (\ref{esol1}), where $A_j$ is satisfied the following transcendental equation
\begin{eqnarray}
\Gamma\left(\frac{\beta_j+\alpha}{1-m_j}+1-\alpha\right)\left(\lambda_j A_j^{m_j}+b_j\right)-\Gamma\left(\frac{\beta_j+\alpha}{1-m_j}+1\right)A_j=0.\label{treq1}
\end{eqnarray}
We suppose that 
\begin{eqnarray}
\nu_j=\frac{\beta_j+m_j\alpha}{1-m_j}. \label{nu}
\end{eqnarray}
Then it is easily verified that the nonlinear equation (\ref{nsfeq2}) has the solution
\begin{eqnarray}
y_j(x)=A_jx^{\frac{\beta_j+\alpha}{1-m_j}}.\label{esol3}
\end{eqnarray}
So we have $\gamma_j=\frac{\beta_j+\alpha}{1-m_j}+1$. It is clear that the exact solution (\ref{esol3}) fulfills not only the boundary conditions (\ref{ebc1}) and satisfies the conditions (\ref{ebc2}). After putting the solution (\ref{esol3}) into the vertex boundary conditions (\ref{bc1})-(\ref{bc2}) we have the following system of transcendental equations with respect to $\lambda_j,\,j=1,2,3$
\begin{eqnarray}
\lambda_1^{\frac{1}{m_1-1}}A_1L_1^{\gamma_1-1}=\lambda_2^{\frac{1}{m_2-1}}A_2L_2^{\gamma_2-1}=\lambda_3^{\frac{1}{m_3-1}}A_3L_3^{\gamma_3-1},\label{streq3}\\
\lambda_1^{\frac{m_1}{m_1-1}}A_1\frac{\Gamma(\gamma_1)}{\Gamma(\gamma_1-\alpha)}L_1^{\gamma_1-\alpha-1}=\lambda_2^{\frac{m_2}{m_2-1}}A_2\frac{\Gamma(\gamma_2)}{\Gamma(\gamma_2-\alpha)}L_2^{\gamma_2-\alpha-1}+\lambda_3^{\frac{m_3}{m_3-1}}A_3\frac{\Gamma(\gamma_3)}{\Gamma(\gamma_3-\alpha)}L_3^{\gamma_3-\alpha-1}.\label{streq4}
\end{eqnarray}
In the particular case when $m_j=m,\,\gamma_j=\gamma,\,L_j=L,\,j=1,2,3$, we have $\lambda_1=\lambda_2+\lambda_3$, therefore the system of transcendental equations (\ref{streq3})-(\ref{streq4}) has a solution.

Now we formulate the following theorem.

{\bf Theorem 2:} The problem (\ref{nsfeq2}) with the boundary conditions (\ref{bc1})-(\ref{ebc1}) has the exact solution in the form (\ref{esol1}).

\section{conclusions}

We considered nonlinear stationary differential equations with fractional-in-space order $1 < \alpha < 2$ on the metric star graph with three finite bonds in this work. We put weight continuity and the generalized Kirchhoff rule at the branched point of the star graph . First we look for the solution on the each bond, after gluing conditions at the vertex and the ends of bonds the exact solutions of nonlinear stationary fractional equations on the star graph were found. These can be extended to the star graphs with any number of bonds.

\end{document}